\documentclass[11pt]{article}

\usepackage{IEEEtrantools}
\usepackage{latexsym}
\usepackage{amsfonts}
\usepackage{amssymb}
\usepackage{amsmath,amsfonts,amsthm}
\usepackage{mathrsfs}
\usepackage{enumerate}
\usepackage{dsfont}

\usepackage{multirow}
\usepackage{tikz}
\usepackage{tikz-cd}
\usetikzlibrary{arrows,decorations.markings}
\usepackage{bookmark}
\usepackage{hyperref}
\hypersetup{
     colorlinks   = true,
     citecolor    = blue
}

\newcommand{\be}{\begin{equation}}
\newcommand{\ee}{\end{equation}}
\newcommand{\bea}{\begin{eqnarray}}
\newcommand{\eea}{\end{eqnarray}}
\newcommand{\bean}{\begin{eqnarray*}}
\newcommand{\eean}{\end{eqnarray*}}
\newcommand{\brray}{\begin{array}}
\newcommand{\erray}{\end{array}}
\newcommand{\biearray}{\begin{IEEEarray}{rCl}}
\newcommand{\eiearray}{\end{IEEEarray}}

\newtheorem{dfn}{Definition}[section]
\newtheorem{thm}[dfn]{Theorem}
\newtheorem{lmma}[dfn]{Lemma}
\newtheorem{ppsn}[dfn]{Proposition}
\newtheorem{crlre}[dfn]{Corollary}
\newtheorem{xmpl}[dfn]{Example}
\newtheorem{rmrk}[dfn]{Remark}

\newcommand{\bdfn}{\begin{dfn}\rm}
\newcommand{\bthm}{\begin{thm}}
\newcommand{\blmma}{\begin{lmma}}
\newcommand{\bppsn}{\begin{ppsn}}
\newcommand{\bcrlre}{\begin{crlre}}
\newcommand{\bxmpl}{\begin{xmpl}}
\newcommand{\brmrk}{\begin{rmrk}\rm}

\newcommand{\edfn}{\end{dfn}}
\newcommand{\ethm}{\end{thm}}
\newcommand{\elmma}{\end{lmma}}
\newcommand{\eppsn}{\end{ppsn}}
\newcommand{\ecrlre}{\end{crlre}}
\newcommand{\exmpl}{\end{xmpl}}
\newcommand{\ermrk}{\end{rmrk}}

\newcommand{\bbz}{\mathbb{Z}}
\newcommand{\bbn}{\mathbb{N}}

\newcommand{\clh}{\mathcal{H}}

\newcommand{\cll}{\mathcal{L}}

\newcommand{\prf}{\noindent{\it Proof\/}: }

\def \qed { \mbox{}\hfill
$\Box$\vspace{1ex}}

\begin{document}          


\author{\sc{Bipul Saurabh}}
\title{  The double covering of   compact quantum group $SO_q(4)$ }
\maketitle


\begin{abstract}
We show that $\mathbb{Z}_2$ is a quantum normal subgroup of $SU_q(2)\times SU_q(2)$ and that the quotient group is isomorphic to $SO_q(4)$, thereby establishing that $SU_q(2)\times SU_q(2)$ is a double cover of $SO_q(4)$.

\end{abstract}
\bigskip

{\bf AMS Subject Classification No.:} {\large 46}L{\large 67},  {\large 20}G{\large 42},  {\large 81}R{\large 50}
\\

\bigskip


\section{Introduction}
In the classical setting, the quotient of the  group $SU(2) \times SU(2)$ by its central normal subgroup $\{I, -I\}$ is isomorphic to the rotation group $SO(4)$. This paper investigates the analogue of this classical result within the framework of compact quantum groups. Specifically, we prove that $\mathbb{Z}_2$ forms a normal quantum subgroup of $SU_q(2) \times SU_q(2)$ and show that the resulting quotient quantum group is isomorphic to the compact quantum group $SO_q(4)$. 

A related construction has been explored in the literature (see \cite{Dij-1994aa},\cite{Lan1998aa}), where $\mathbb{Z}_2$ was shown to be a normal quantum subgroup of $SU_{\sqrt{q}}(2)$ by mapping the standard generators $a$ and $b$ to the generating unitary and zero of $C(\mathbb{Z}_2)$, respectively. The corresponding quotient is isomorphic to $SO_q(3)$, thereby establishing $SU_{\sqrt{q}}(2)$ as the quantum double covering of $SO_q(3)$. However, to the best of our knowledge, the quantum double covering of $SO_q(4)$ has not yet been addressed. This work aims to fill this gap in the literature.

Throughout the paper,  $q$ will denote a real number in $(0,1)$.   If $M$ is a subset of a $C^*$-algebra $A$, then $\langle M\rangle$ denotes the closure of the linear span of $M$. Whenever we consider the tensor product $A\otimes B$ of two $C^*$-algebras $A$ and $B$, we mean the spatial tensor product $C^*$-algebra. In the paper,  one of the factors is always nuclear, and hence all $C^*$-algebra tensor product norms coincide.
    \section{Preliminaries}
 
Let $\mathbb{G}=(C(\mathbb{G}),\Delta_{\mathbb{G}})$ and
$\mathbb{H}=(C(\mathbb{H}),\Delta_{\mathbb{H}})$ be compact quantum
groups. We say that $\mathbb{H}$ is a \emph{quantum subgroup} of
$\mathbb{G}$ if there exists a surjective unital $*$-homomorphism
\[
\pi:C(\mathbb{G})\longrightarrow C(\mathbb{H})
\]
such that
\[
\Delta_{\mathbb{H}}\circ\pi
=
(\pi\otimes\pi)\circ\Delta_{\mathbb{G}}.
\]
In other words, $\pi$ is a morphism of compact quantum groups. 
Further, we say that $\mathbb{H}$ is a normal quantum subgroup of
$\mathbb{G}$ if
\[
(\operatorname{id}\otimes h_{\mathbb{H}}\circ\pi)
\Delta_{\mathbb{G}}(a)
=
(h_{\mathbb{H}}\circ\pi\otimes\operatorname{id})
\Delta_{\mathbb{G}}(a),
\qquad a\in C(\mathbb{G}),
\]
where $h_{\mathbb{H}}$ is the Haar state on $C(\mathbb{H})$.
Suppose that $\mathbb{H}$ is a normal quantum subgroup of $\mathbb{G}$, with quotient morphism
\[
\pi:C(\mathbb{G})\longrightarrow C(\mathbb{H}).
\]
The corresponding quotient compact quantum group, denoted by
$\mathbb{G}/\mathbb{H}$, can be realized as the compact quantum group
whose $C^*$-algebra is
\[
C(\mathbb{G}/\mathbb{H})
=
\left\{
a\in C(\mathbb{G}) :
(\operatorname{id}\otimes h_{\mathbb{H}}\circ\pi)
\Delta_{\mathbb{G}}(a)=a
\right\}.
\]
The comultiplication on $C(\mathbb{G})$ restricts to a comultiplication on
$C(\mathbb{G}/\mathbb{H})$, making it a compact quantum group. Its Haar
state $h_{\mathbb{G}/\mathbb{H}}$ is characterized by
\[
(h_{\mathbb{G}/\mathbb{H}}\otimes\operatorname{id})
\Delta_{\mathbb{G}/\mathbb{H}}(a)
=
h_{\mathbb{G}/\mathbb{H}}(a)1
=
(\operatorname{id}\otimes h_{\mathbb{G}/\mathbb{H}})
\Delta_{\mathbb{G}/\mathbb{H}}(a),
\qquad
a\in C(\mathbb{G}/\mathbb{H}).
\]
Let $\operatorname{Pol}(\mathbb{G}) \subset C(\mathbb{G})$ be  the linear span of the matrix coefficients of all finite-dimensional unitary corepresentations of $\mathbb{G}$. 
By Woronowicz's Peter--Weyl theorem, $\operatorname{Pol}(\mathbb{G})$ forms a dense, canonical Hopf $*$-subalgebra of $C(\mathbb{G})$ under the restriction of the coproduct $\Delta$, with the coundit $\epsilon$ and antipode $S$ uniquely well-defined on it.

\bppsn \label{Hopf}
Let $\mathbb{G}$ and $\mathbb{H}$ be compact quantum groups with
faithful Haar states $h_{\mathbb{G}}$ and $h_{\mathbb{H}}$, respectively.
Suppose that
\[
\phi:C(\mathbb{G})\longrightarrow C(\mathbb{H})
\]
is a surjective unital $*$-homomorphism satisfying
\[
\Delta_{\mathbb{H}}\circ\phi
=
(\phi\otimes\phi)\circ\Delta_{\mathbb{G}}.
\]
If the restriction
\[
\phi_0:=\phi|_{\operatorname{Pol}(\mathbb{G})}:
\operatorname{Pol}(\mathbb{G})
\longrightarrow
\operatorname{Pol}(\mathbb{H})
\]
is an isomorphism of Hopf $*$-algebras, then $\phi$ is a
$*$-isomorphism.
\eppsn

\prf
Since $\phi_0$ is a Hopf $*$-algebra isomorphism, it preserves the coproduct, implying that the functional $h_{\mathbb{H}}\circ\phi_0$ is a Haar functional on $\operatorname{Pol}(\mathbb{G})$. By the uniqueness of the Haar functional, $h_{\mathbb{H}}\circ\phi_0=h_{\mathbb{G}}$ on $\operatorname{Pol}(\mathbb{G})$.

Consequently, $\phi_0$ preserves the $L^2$-inner product. Since the Haar states are assumed to be faithful, the canonical GNS representations are faithful, allowing us to identify $C(\mathbb{G})$ and $C(\mathbb{H})$ with their respective images acting on the Hilbert spaces $L^2(\mathbb{G})$ and $L^2(\mathbb{H})$. Because $\phi_0$ is a bijection between the dense subspaces $\operatorname{Pol}(\mathbb{G})$ and $\operatorname{Pol}(\mathbb{H})$ that preserves the inner product, it extends uniquely to a unitary operator $U: L^2(\mathbb{G}) \to L^2(\mathbb{H})$ satisfying $\phi_0(a) = U a U^*$ for all $a \in \operatorname{Pol}(\mathbb{G})$. It immediately follows that $\phi_0$ preserves the operator norm:
\[
\|\phi_0(a)\|_{C(\mathbb{H})} = \|U a U^*\| = \|a\|_{C(\mathbb{G})}
\]
for all $a \in \operatorname{Pol}(\mathbb{G})$. By the density of $\operatorname{Pol}(\mathbb{G})$ in $C(\mathbb{G})$ and the continuity of $\phi$, the unique isometric extension of $\phi_0$ to $C(\mathbb{G})$ must coincide with $\phi$. Since $\phi$ is an isometry, it is injective.
\qed

\section{$C(SO_q(4)$)}	In this subsection, we provide a brief overview of the Hopf algebra structure associated with the compact quantum group $SO_q(4)$ as introduced in \cite[Section 9.3.3]{KliSch-1997aa}. We begin by introducing key notations. For $1\leq i,j,m,n\leq 4$ we define
\begin{IEEEeqnarray}{rCl} 
	&i'=& 5-i,
	~~\rho_{i}  =  2-i \text{ if } i < i',
	~~\rho_{i'} = -\rho_i \text{ if } i \leq i',
	~~C_j^i =  \delta_{ij'} q^{-\rho_i},\nonumber\\
	&R_{mn}^{ij}=&\begin{cases}
		(q-q^{-1})(\delta_{jm}\delta_{in}- C_{i}^{j}C_{n}^{m}) & \mbox{ if } i>m, \cr
		q^{\delta_{ij}-\delta_{ij^{'}}}\delta_{im}\delta_{jn}  & \mbox{ if } i\leq m, \cr
	\end{cases},\nonumber
\end{IEEEeqnarray}
where $\delta_{ij}$ represents the Kronecker delta function. Let $A(R)$ be the unital associative algebra generated by $v_j^i$, $i,j=1,,\ldots,4$, subject to the following relations:
\begin{IEEEeqnarray}{rCl} \label{Relation 2.1}
	\sum_{k,l=1}^{4}R_{kl}^{ji}v_{s}^{k}v_{t}^{l}-R_{st}^{lk}v_{k}^{i}v_{l}^{j}=0,
	\quad i,j,s,t=1,\ldots , 4.
\end{IEEEeqnarray}
The matrices $(\!(v_j^i)\!)$ and $(\!(C_j^i)\!)$ are denoted as $V$ and $C$, respectively. Define $J$ as the two-sided ideal of $A(R)$ generated by the entries of the matrices $VCV^t C^{-1}-I$ and $CV^t C^{-1}V-I$. Let $\mathcal{O}\left(O_q(4)\right)$ denotes the quotient algebra $A(R)/J$. The Hopf $*$-algebra structure on $\mathcal{O}(O_{q}(4)$ comes from the following maps.
\begin{IEEEeqnarray*}{lll}
	&\bullet \, \text{Comultiplication}
	:\Delta(v_{l}^{k}) = \sum_{i=1}^{4}v_{i}^{k} \otimes v_{l}^{i}, \quad \quad & \bullet \, \text{Counit}:
	\epsilon(v_{l}^{k}) = \delta_{kl}, \\
	& \bullet\, \text{Antipode}:
	S(v_l^k) = q^{\rho_k-\rho_l}v_{l'}^{k'}, &\bullet \, \text{Involution}:
	(v_{l}^{k})^* = q^{\rho_k - \rho_l}v_{l'}^{k'}.
\end{IEEEeqnarray*}

Let $\mathcal{D}_q$ be the quantum determinant of the matrix $V$ for the quantum group $O_q(4)$ \cite[Chapter 9, Definition 10]{KliSch-1997aa}. Denote the quotient $\mathcal{O}(O_q(4))/\langle \mathcal{D}_q - 1\rangle$ of $\mathcal{O}(O_q(4))$ by the two-sided ideal $\langle\mathcal{D}_q - 1\rangle$ as $\mathcal{O}(SO_q(4))$. The Hopf $*$-algebra structure on $\mathcal{O}(SO_q(4))$ is induced from $\mathcal{O}(O_q(4))$. In $\mathcal{O}\left(SO_q(4)\right)$, the relation $V^*=CV^t C^{-1}$ leads to the following:
\begin{equation}\label{Relation 2.2}
	VV^* = V^*V = I.
\end{equation}
The algebra $\mathcal{O}\left(SO_q(4)\right)$ becomes a normed $*$-algebra with the norm defined by
$$\|a\|=\sup\left\{\|\pi(a)\|: \pi \text{ is a representation of } \mathcal{O}\left(SO_q(4)\right)\right\} \text{ for } a\in \mathcal{O}(SO_q(4)).$$

We denote the completion of $\mathcal{O}\left(SO_q(4)\right)$ with respect to this norm as $C\left(SO_q(4)\right)$. The pair $\left(C\left(SO_q(4)\right),\Delta\right)$ forms a compact quantum group, known as the $q$-deformation of the group $SO(4)$. For further details, we refer the reader to \cite[Section 9.3.3]{KliSch-1997aa}.
\subsection{Relations:} 
Here we describe all the defining relations of $C(SO_q(4))$ explicitly using equation $(\ref{Relation 2.1},\ref{Relation 2.2})$ and the formula for the quantum determinant $\mathcal{D}_q$ given in  \cite{Fio-1994aa}.  Set $z_i=v^4_{5-i}$ and $y_i=v^3_{5-i}$ for $1\leq i \leq 4$. 
\subsection*{(1). Relations within the $z$-Sector}

\subsubsection*{Relations among $z_i$}
\begin{align*}
	z_2 z_1 &= q z_1 z_2,  & z_3 z_2 &= z_2 z_3, \\
	z_3 z_1 &= q z_1 z_3,  & z_4 z_2 &= q z_2 z_4, \\
	z_4 z_1 &= q^2 z_1 z_4, & z_4 z_3 &= q z_3 z_4.
\end{align*}

\subsubsection*{Relations between $z_i$ and $z_j^*$}
\begin{align*}
	z_4 z_4^* &= z_4^* z_4 \\
	z_3 z_3^* &= z_3^* z_3 + (q - q^{-1}) z_4^* z_4 \\
	z_2 z_2^* &= z_2^* z_2 + (q - q^{-1}) z_4^* z_4 \\
	z_1 z_1^* &= z_1^* z_1 + (q - q^{-1}) z_2^* z_2 + (q - q^{-1}) z_3^* z_3 + (q - q^{-1})^2 z_4^* z_4 \\[0.5em]
	z_4 z_3^* &= q^{-1} z_3^* z_4,        & z_3 z_4^* &= q z_4^* z_3, \\
	z_4 z_2^* &= q^{-1} z_2^* z_4,        & z_2 z_4^* &= q z_4^* z_2, \\
	z_4 z_1^* &= q^{-2} z_1^* z_4,        & z_1 z_4^* &= q^2 z_4^* z_1, \\[0.5em]
	\intertext{\noindent\small\itshape Cross-commutators with drift terms:}
	z_3 z_2^* &= z_2^* z_3 + (q - q^{-1}) z_4^* z_1, \\
	z_2 z_3^* &= z_3^* z_2 + (q - q^{-1}) z_4^* z_1, \\[0.5em]
	z_3 z_1^* &= q^{-1} z_1^* z_3,        & z_1 z_3^* &= q z_3^* z_1, \\
	z_2 z_1^* &= q^{-1} z_1^* z_2,        & z_1 z_2^* &= q z_2^* z_1.
\end{align*}

\subsection*{(2). Relations within the $y$-Sector}

\subsubsection*{Relations among $y_i$}
\begin{align*}
	y_2 y_1 &= q y_1 y_2,  & y_4 y_2 &= q y_2 y_4, \\
	y_3 y_1 &= q y_1 y_3,  & y_4 y_3 &= q y_3 y_4, \\
	y_3 y_2 &= y_2 y_3,    & y_1 y_4 &= q^2 y_4 y_1 + (q^2 - 1) y_3 y_2.
\end{align*}

\subsubsection*{Relations between $y_i$ and $y_j^*$}
\begin{align*}
	y_3 y_3^* &= y_3^* y_3 \\
	y_4 y_4^* &= y_4^* y_4 + (1 - q^2) y_3^* y_3 \\
	y_1 y_1^* &= y_1^* y_1 + (1 - q^2) y_3^* y_3 \\
	y_2 y_2^* &= q^2 y_2^* y_2 + (1 - q^2) y_4^* y_4 + (1 - q^2) y_1^* y_1 + (q^2 - 2 + q^{-2}) y_3^* y_3 \\[0.5em]
	y_3 y_4^* &= q^{-1} y_4^* y_3, & y_4 y_3^* &= q y_3^* y_4, \\
	y_3 y_1^* &= q^{-1} y_1^* y_3, & y_1 y_3^* &= q y_3^* y_1, \\
	y_3 y_2^* &= q^{-2} y_2^* y_3, & y_2 y_3^* &= y_3^* y_2, \\[0.5em]
	y_4 y_1^* &= y_1^* y_4,        & y_1 y_4^* &= y_4^* y_1, \\
	y_4 y_2^* &= q^{-1} y_2^* y_4, & y_2 y_4^* &= q y_4^* y_4 - (q - q^{-1}) y_3^* y_1, \\
	y_1 y_2^* &= q^{-1} y_2^* y_1, & y_2 y_1^* &= q y_1^* y_2 - (q - q^{-1}) y_3^* y_4.
\end{align*}

\subsection*{(3). Cross-Relations }

\subsubsection*{Relations between $y_i$ and $z_j$}
\begin{align*}
	y_4 z_4 &= q z_4 y_4,                  & y_4 z_2 &= z_2 y_4, \\
	y_4 z_3 &= z_3 y_4,                  & y_4 z_1 &= q^{-1} z_1 y_4, \\
	y_2 z_2 &= q z_2 y_2,                  & y_2 z_1 &= z_1 y_2, \\
	y_3 z_3 &= q z_3 y_3,                  & y_3 z_1 &= z_1 y_3, \\[0.5em]
	y_2 z_4 &= (q - q^{-1}) z_2 y_4 + z_4 y_2, \\
	y_2 z_3 &= (q^{-2} - 1) z_1 y_4 + q^{-1} z_3 y_2, \\
	y_3 z_4 &= (q - q^{-1}) z_3 y_4 + z_4 y_3, \\
	y_3 z_2 &= (q^{-2} - 1) z_1 y_4 + q^{-1} z_2 y_3, \\
	y_1 z_1 &= q z_1 y_1, \\
	y_1 z_2 &= (q - q^{-1}) z_1 y_2 + z_2 y_1, \\
	y_1 z_3 &= (q - q^{-1}) z_1 y_3 + z_3 y_1, \\
	y_1 z_4 &= (q - 2q^{-1} + q^{-3}) z_1 y_4 + (q^{-2} - 1) z_2 y_3 \\
	&\quad + (q^{-2} - 1) z_3 y_2 + q^{-1} z_4 y_1.
\end{align*}

\subsubsection*{Relations between $z_i^*$ and $y_j$}
\begin{align*}
	z_4^* y_4 &= q y_4 z_4^*,                  & z_4^* y_2 &= y_2 z_4^*, \\
	z_4^* y_3 &= y_3 z_4^*,                  & z_4^* y_1 &= q y_1 z_4^*, \\
	z_1^* y_4 &= q^{-1} y_4 z_1^*,            & z_2^* y_4 &= y_4 z_2^*, \\
	z_3^* y_4 &= y_4 z_3^*,                  & z_3^* y_2 &= y_2 z_3^*, \\
	z_2^* y_3 &= y_3 z_2^*,                  & & \\[0.5em]
	z_1^* y_3 &= y_3 z_1^* - (q - q^{-1}) y_4 z_2^*, \\
	z_1^* y_2 &= y_2 z_1^* - (q - q^{-1}) y_4 z_3^*, \\
	z_2^* y_1 &= q z_1 y_2^*, \\
	z_3^* y_1 &= q z_1 y_3^*, \\[0.5em]
	z_2^* y_2 &= q^{-1} y_2 z_2^* - (q - q^{-1}) y_4 z_4^*, \\
	z_3^* y_3 &= q^{-1} y_3 z_3^* - (q - q^{-1}) y_4 z_4^*, \\[0.5em]
	z_1^* y_1 &= q y_1 z_1^* - (q^2 - 2 + q^{-2}) y_4 z_4^* 
 - (q^{-1} - q) y_3 z_3^* - (q^{-1} - q) y_2 z_2^*.
\end{align*}
\subsection*{(4). Orthogonality Relations from Equation (2.2)}

\begin{align*}
	\sum_{m=1}^{4} z_m z_m^* &= 1, \qquad\qquad \sum_{m=1}^{4} y_m y_m^* = 1, \\
	\sum_{m=1}^{4} z_m y_m^* &= 0, \qquad\qquad \sum_{m=1}^{4} y_m z_m^* = 0.
\end{align*}

\begin{align*}
	z_i^* z_j + y_i^* y_j + q^{\rho_{5-i} + \rho_{5-j}} \left( y_{5-i} y_{5-j}^* + q^2 z_{5-i} z_{5-j}^* \right) = \delta_{ij}.
\end{align*}

\subsection*{(5). The quantum determinant  $\mathcal{D}_q$}
\begin{align*}
	\mathcal{D}_q := \;& z_4^* \Bigl[ y_3^* y_3 z_4 - q^{-1} y_3^* y_4 z_3 - q^{-1} y_2^* y_2 z_4 
 + y_2^* y_4 z_2 - q^{-1} y_1^* y_2 z_3 + y_1^* y_3 z_2 \Bigr] \\
	+ \;& z_3^* \Bigl[ -q^{-1} y_4^* y_1 z_3 + y_4^* y_4 z_3 + q y_4^* y_2 z_1 
 - y_1^* y_4 z_2 + q y_2^* y_2 z_3 - q^2 y_2^* y_3 z_3 \Bigr] \\
	+ \;& z_2^* \Bigl[ -q^{-1} y_4^* y_1 z_2 + q^{-2} y_4^* y_4 z_2 + q y_4^* y_3 z_1 
- y_1^* y_4 z_3 - q^2 y_3^* y_3 z_2 + q^3 y_3^* y_2 z_2 \Bigr] \\
	+ \;& z_1^* \Bigl[ q y_4^* y_2 z_3 - y_4^* y_3 z_2 - q^2 y_3^* y_1 z_3 
 + q^3 y_3^* y_3 z_1 + q^3 y_2^* y_1 z_2 - q^4 y_2^* y_2 z_1 \Bigr] = 1.
\end{align*}
We refer to the set of relations mentioned above as $\mathcal{R}_4$. Thus, the $C^*$-algebra $C(SO_q(4))$ is the universal $C^*$-algebra generated by $\{z_i, y_i:1\leq i\leq 4\}$ subject to the relations $\mathcal{R}_4$.

\subsection{Faithful representation} Let $\clh= \ell^2(\bbz)\otimes  \ell^2(\bbz)  \otimes   \ell^2(\bbn)\otimes  \ell^2(\bbn)$. Let $\alpha:=S^*\sqrt{1-q^{2N}}$ and $\beta:=q^{N}$. Define the representation $\Psi: C\left(SO_q(4)\right) \rightarrow \cll(\clh)$  given by 
\begin{IEEEeqnarray*}{rClrCl}
	\Psi(z_1)&=& t_1 \otimes 1  \otimes   \alpha^* \otimes  \alpha^*, \qquad & 
	\Psi(y_1)&=& 1  \otimes t_2 \otimes  q^{N+1}\otimes \alpha^*,\\
	\Psi(z_2)&=& t_1 \otimes 1  \otimes  q^N\otimes \alpha^*,	\qquad &
	\Psi(y_2)&=&  1  \otimes t_2   \otimes \alpha \otimes \alpha^*\\
	\Psi(z_3)&=&  t_1 \otimes 1  \otimes  \alpha^*\otimes q^N, \qquad &
	\Psi(y_3)&=& 1  \otimes t_2 \otimes q^N\otimes q^{N+1},\\
	\Psi(z_4)&=& t_1 \otimes 1  \otimes q^N\otimes q^N,  \qquad & 
	\Psi(y_4)&=& 1  \otimes t_2 \otimes \alpha \otimes q^N.
\end{IEEEeqnarray*}
It follows from \cite{BhuBisSau-2024aa} that $\Psi$  is a faithful representation. Henceforth, we identify $z_{i}$ and $y_{i}$ with their images under $\Psi$.

	\section{The quotient group  $SU_q(2)\times SU_q(2)/\bbz_2$}
In this section, we prove out main result. We start with a definition. 
\bdfn 
Let $A_q$ denote the $C^*$-algebra of continuous functions on $SU_q(2)$, which is the universal $C^*$-algebra generated by two elements $a$ and $b$ subject to the relations 
\[
a^*a + b^*b = 1= aa^* + q^2bb^*,
\]
\[
ab= qba, \quad ab^* =qb^*a, \quad \text{and} \quad b^*b = bb^*.
\]
The group structure is given by the coproduct $\Delta : A_q \to A_q \otimes A_q$ defined by
\begin{align*}
	\Delta(a) &= a \otimes a - qb^* \otimes b, \\
	\Delta(b) &= b \otimes a + a^* \otimes b.
\end{align*}

The faithful Soibelman representation of $A_q$ is given by 
\[
\pi : A_q \longrightarrow \mathcal{L}\bigl(\ell^2(\mathbb{N}_0) \otimes \ell^2(\bbz)\bigr);\quad  \pi(a)=S\sqrt{1-q^{2N}} \otimes 1, \mbox{ and } \pi(b)=q^N \otimes S.
\]
The coordinate algebra of the quantum group $SU_q(2)$, denoted by $\mathcal{O}(SU_q(2))$, is the $\star$-subalgebra of $A_q$ generated by $a$ and $b$. The fundamental generating unitary is given by 
\[
U:=\begin{pmatrix}
a & -qb^*\\
b & a^*
\end{pmatrix}.
\]
The  quantum group $SU_q(2) \times SU_q(2)$ is defined such that its $C^*$-algebra of continuous functions is the minimal tensor product $A_q \otimes A_q$, and the   coproduct $$\Delta_{\times} : A_q \otimes A_q \to (A_q \otimes A_q) \otimes (A_q \otimes A_q)$$ is given by $\Delta_{\times} = (\text{id} \otimes \Sigma \otimes \text{id}) \circ (\Delta \otimes \Delta)$, where $\Sigma : A_q \otimes A_q \to A_q \otimes A_q$ is the tensor flip map $\Sigma(a \otimes b) = b \otimes a$.
\edfn
It is well known that the equivalence classes of irreducible unitary corepresentations of $SU_q(2)$ are parametrized by half-integers. Let 
$$\operatorname{Irr}(SU_q(2)) = \{u^{(j)} : j \in \frac{1}{2}\mathbb{N}_0 = \{0, \frac{1}{2}, 1, \frac{3}{2}, \dots\}\},$$
 where each $u^{(j)}$ acts on a $(2j+1)$-dimensional Hilbert space $V_j$. 
Let $u^{(j)}$ and $u^{(k)}$ be two such irreducible corepresentations. The \textbf{external tensor product} $u^{(j)}\boxtimes u^{(k)}$ is a corepresentation of  $SU_q(2)\times SU_q(2)$ acting on  $V_j\otimes V_k$. In terms of the spatial leg-numbering notation within the algebra $\mathcal{B}(V_j \otimes V_k) \otimes C(SU_q(2) \times SU_q(2)) \cong \mathcal{B}(V_j) \otimes \mathcal{B}(V_k) \otimes C(SU_q(2)) \otimes C(SU_q(2))$, it is defined by
\[
u^{(j)}\boxtimes u^{(k)} = u^{(j)}_{13} u^{(k)}_{24}.
\]
More explicitly, if we fix orthonormal bases for $V_j$ and $V_k$ so that the corepresentations are written as matrices of coefficients $u^{(j)} = (u^{(j)}_{mn})_{m,n=-j}^j$ and $u^{(k)} = (u^{(k)}_{pq})_{p,q=-k}^k$, the matrix coefficients of the external tensor product on $V_j \otimes V_k$ are given by 
\[
(u^{(j)}\boxtimes u^{(k)})_{(m,p),(n,q)} = u^{(j)}_{mn} \otimes u^{(k)}_{pq} \,\in\, C(SU_q(2)) \otimes C(SU_q(2)) \cong C(SU_q(2) \times SU_q(2)).
\]
Also, we have 
\[
\dim(u^{(j)}\boxtimes u^{(k)})
=(2j+1)(2k+1).
\]
Every irreducible corepresentation of $SU_q(2)\times SU_q(2)$
is obtained in this way. Hence
\[
\operatorname{Irr}(SU_q(2)\times SU_q(2))
=
\{u^{(j)}\boxtimes u^{(k)}:
j,k=0,\tfrac12,1,\ldots\}.
\]
Thus,  coordinate algebra $\mathcal{O}(SU_q(2)\times SU_q(2))$ is the algebraic tensor product $\mathcal{O}(SU_q(2)) \otimes \mathcal{O}(SU_q(2))$. 
\bdfn
 The $C^*$-algebra of the compact quantum group $(C(\mathbb{Z}_2), \Delta_{\mathbb{Z}_2})$ is generated by a single self-adjoint unitary function $z$ satisfying $z([0]) = 1$, $z([1])=-1$,  and $z^2 = 1$. The coproduct is given by $\Delta_{\mathbb{Z}_2}(z) = z \otimes z$, while the counit and antipode are  defined by $\varepsilon(z) = 1$ and $S(z) = z$, respectively.
 \edfn
Let's first  briefly recall from \cite{Dij-1994aa} the quotient group $SU_q(2)/\bbz_2$. 
For that, define a  surjective $C^*$-homomorphism 
$$ \phi: A_q\rightarrow C(\bbz_2); \, \, U\mapsto \begin{pmatrix}
	z & 0 \\
	0 & z
\end{pmatrix}.$$ 
This establish $\bbz_2$ a normal quantum subgroup of $SU_q(2)$.  We get an action of $\bbz_2$ on $A_q$ as follows:
$$\varphi: A_q\rightarrow A_q \otimes C(\bbz_2); \,\, \varphi(x)=(1\otimes \varphi)\Delta(x), \mbox{ for }x \in A_q.$$
Equivalently we have a group homomorphism  
$$\gamma: \bbz_2 \rightarrow \mbox{ Aut}(A_q); \quad  \gamma(g)(x)= (1\otimes \mbox{ev}_g)\varphi(x).$$
The automorphism associated to the generating element $[1]$ is given by  
$$\gamma([1])(\alpha)=-\alpha, \mbox{ and } \gamma([1])(\beta)=-\beta.$$
In such a situation, one defines the quotient space $SU_q(2)/\bbz_2$ as follows:
$$C(SU_q(2)/\bbz_2):=\{x \in A_q: \phi(x)=x \otimes 1\}=\{x \in A_q: \gamma([1])(x)=x\}=A_q^0.$$
It is proved in \cite{Dij-1994aa} that 
$$A_q^0=CLS\langle C^n_{ij}: -n \leq i,j \leq n, n \in \bbn_0,  \rangle$$
establishing that $SU_q(2)/\bbz_2$ is isomorphic to $SO_q(3)$.

Let $m: C(\mathbb{Z}_2) \otimes C(\mathbb{Z}_2) \rightarrow C(\mathbb{Z}_2)$ be the standard pointwise multiplication map. Since $C(\mathbb{Z}_2)$ is commutative, $m$ is a well-defined $*$-homomorphism.
Define a  surjective $C^*$-homomorphism 
$$ \phi_{\times}: A_q \otimes A_q\xrightarrow{\phi\otimes \phi} C(\bbz_2)\otimes  C(\bbz_2) \xrightarrow{m} C(\bbz_2); \, \, \quad  \phi_{\times}(x\otimes y)=\phi(x)\phi(y) \mbox{ for } x,y \in A_q \otimes A_q.$$ 
In particular, we have 
$$ \phi_{\times}(U\otimes 1) =\phi_{\times}( 1 \otimes U)= \begin{pmatrix}
	z & 0 \\
	0 & z
\end{pmatrix}. $$

\bppsn \label{normal}  The group $\mathbb{Z}_2$ is a normal quantum subgroup of $SU_q(2)\times SU_q(2)$.  
\eppsn
\begin{proof} 
	To establish that $\mathbb{Z}_2$ is a quantum subgroup, it suffices to prove that the relation 
	\[ 
	\Delta_{\mathbb{Z}_2}\bigl(\phi_{\times}(x)\bigr) = (\phi_{\times} \otimes \phi_{\times})\bigl(\Delta_{\times}(x)\bigr) 
	\] 
	holds for all $x \in \{u_{ij} \otimes 1, 1\otimes u_{ij}: -\frac{1}{2}\leq i,j \leq \frac{1}{2}\}$, as these elements generate the dense $*$-subalgebra $\operatorname{Pol}(SU_q(2) \times SU_q(2))$. 
For that, we have
	\begin{IEEEeqnarray*}{rCl} 
		(\phi_{\times} \otimes \phi_{\times})\bigl(\Delta_{\times}(u_{ij} \otimes 1)\bigr) &=& (\phi_{\times} \otimes \phi_{\times})\left( \sum_{k=-\frac{1}{2}}^{\frac{1}{2}} (u_{ik} \otimes 1) \otimes (u_{kj} \otimes 1) \right) \\ 
		&=& \sum_{k=-\frac{1}{2}}^{\frac{1}{2}} \phi_{\times}(u_{ik} \otimes 1) \otimes \phi_{\times}(u_{kj} \otimes 1) \\ 
		&=& \sum_{k=-\frac{1}{2}}^{\frac{1}{2}} (z \delta_{ik}) \otimes (z \delta_{kj}) \\ 
		&=& z \otimes z \cdot \delta_{ij}. 
	\end{IEEEeqnarray*} 
	On the other hand, we have 
	\begin{IEEEeqnarray*}{rCl} 
		\Delta_{\mathbb{Z}_2}\bigl(\phi_{\times}(u_{ij} \otimes 1)\bigr) &=& \Delta_{\mathbb{Z}_2}(z \delta_{ij}) 
		= \delta_{ij} \Delta_{\mathbb{Z}_2}(z)  
		= z \otimes z \cdot \delta_{ij}. 
	\end{IEEEeqnarray*} 
	To prove the rest, take $\ell_1,\ell_2 \in \frac{1}{2}\mathbb{N}_0$. Then we have 
	\begin{IEEEeqnarray*}{rCl} 
		\phi_{\times}\left((u^{(\ell_1)}\boxtimes u^{(\ell_2)})_{(a,c),(b,d)} \right) &=& \phi_{\times}\left(u^{(\ell_1)}_{ab}\otimes u^{(\ell_2)}_{cd}\right) \\ 
		&=& \phi (u^{(\ell_1)}_{ab})\phi(u^{(\ell_2)}_{cd}) \\ 
		&=& z^{2\ell_1 + 2\ell_2} \cdot I_{(2\ell_1+1)(2\ell_2+1)}. 
	\end{IEEEeqnarray*} 
	Therefore, we obtain
	\begin{IEEEeqnarray}{rCl} \label{eqmul}
		(h_{\mathbb{Z}_2} \otimes \text{id})\left(\phi_{\times}\left(u^{(\ell_1)}\boxtimes u^{(\ell_2)}\right)\right) & = & h_{\mathbb{Z}_2}\left(z^{2\ell_1+2\ell_2}\right) \cdot I_{(2\ell_1+1)(2\ell_2+1)} \nonumber \\ 
		& = & \begin{cases} I_{(2\ell_1+1)(2\ell_2+1)} & \text{if } \ell_1 + \ell_2 \in \mathbb{N}, \\ 0 & \text{if } \ell_1 + \ell_2 \in \mathbb{N}+\frac{1}{2}. \end{cases} 
	\end{IEEEeqnarray} 
Applying this, the claim follows from \cite{Wan-1995aa}.
\end{proof}

Now we get an action of $\bbz_2$ on $A_q$ as follows:
$$\varphi: A_q\otimes A_q \rightarrow ( A_q \otimes A_q )\otimes C(\bbz_2); \,\, \varphi_{\times}(x)=(1\otimes \phi_{\times})\Delta_{\times}(x), \mbox{ for }x \in A_q \otimes A_q.$$
Equivalently we have a group homomorphism  
$$\gamma_{\times}: \bbz_2 \rightarrow \mbox{ Aut}(A_q \otimes A_q); \quad  \gamma_{\times}(g)(x)= (1\otimes \mbox{ev}_g)\varphi_{\times}(x), \mbox{ for }x \in A_q \otimes A_q.$$
The automorphism associated to the generating element $[1]$ is given on the simple tensors by 
$$\gamma_{\times}([1])(x\otimes y )=\gamma([1])(x)\otimes \gamma([1])(y).$$

\blmma \label{lmma1} Let $C^{(\ell_1,\ell_2)}_{(ij)\otimes (kl)}$ be the $(ij)\otimes (kl)$-th matrix element  of the irreducible corepresentation $u^{(\ell_1)}\boxtimes u^{(\ell_2)}$. One has 
$$ C(SU_q(2)\otimes SU_q(2)/\mathbb{Z}_2)=\langle C^{(\ell_1,\ell_2)}_{ij}: \ell_1,\ell_2 \in \frac{1}{2}\bbn_0, \,  \ell_1+\ell_2\in \bbn_0, \,  \ell_1 \leq i \leq \ell_1, \, -\ell_2 \leq j \leq \ell_2 \rangle. $$
\elmma 
\prf By Theorem 1.5 in  \cite{Pod-1995aa}), we get,
\begin{IEEEeqnarray}{rCl}
C(SU_q(2)\otimes SU_q(2)/\mathbb{Z}_2)= \overline{\oplus_{(\ell_1,\ell_2) \in \frac{1}{2}\bbn_0}\oplus_{i \in M_{(\ell_1,\ell_2)} }W_{(\ell_1,\ell_2),i}}, \nonumber 
\end{IEEEeqnarray}
where $(\ell_1,\ell_2)$ represents the irreducible co-representation $u^{(\ell_1)}\boxtimes u^{(\ell_2)}$ of $SU_q(2)\otimes SU_q(2)$, $W_{(\ell_1,\ell_2), i}$ denotes the spectral subspace corresponding to 
$u^{(\ell_1)}\boxtimes u^{(\ell_2)}$ for all $i \in M_{(\ell_1,\ell_2)}$ and $M_{(\ell_1,\ell_2)}$ is the multiplicity of $u^{(\ell_1)}\boxtimes u^{(\ell_2)}$.  It follows from equation (\ref{eqmul}) that 
$$M_{(\ell_1,\ell_2)}=
	\begin{cases} (2\ell_1+1)(2\ell_2+1) & \text{if } \ell_1 + \ell_2 \in \mathbb{N}_0, \\ 
		0 & \text{if } \ell_1 + \ell_2 \in \mathbb{N}_0+\frac{1}{2}. 
	\end{cases}
$$ 
This proves the claim. 
\qed 

\bppsn \label{p1} The $C^*$-algebra $C(SU_q(2)\otimes SU_q(2)/\mathbb{Z}_2)$ is generated by the elements of the unitary corepresentation matrix $U\boxtimes U$. 
\eppsn
\prf
	First, note from Lemma \ref{lmma1} that all matrix elements of the unitary corepresentation matrix $U\boxtimes U$ are in $C(SU_q(2)\otimes SU_q(2)/\mathbb{Z}_2)$. Let $\mathcal{B}$ be the $C^*$-subalgebra of $C(SU_q(2)\otimes SU_q(2)/\mathbb{Z}_2)$ generated by the matrix elements of $U\boxtimes U$. Note that 
	$$1\otimes U^2 = (U U^* \otimes U^2) = (U \boxtimes U) (U^* \boxtimes U).$$
	Therefore, the matrix elements of $1\otimes U^2$ lie in the $C^*$-algebra $\mathcal{B}$. We have two possible cases for the spin labels $(\ell_1, \ell_2)$:
	\begin{itemize}
		\item \textbf{Case 1 ($\ell_1, \ell_2 \in \mathbb{Z}$):} In this case, assuming $n \geq m$, we have 
		\begin{equation}
			(U \otimes 1)^{2m}(1 \otimes U)^{2n} = (U \boxtimes U)^{2m} (1 \otimes U^2)^{n-m}.
		\end{equation}
		
		\item \textbf{Case 2 ($\ell_1, \ell_2 \in \mathbb{Z} + \frac{1}{2}$):} In this case, we have 
		\begin{equation}
			U^{2m+1} \otimes U^{2n+1} = (U \boxtimes U) \cdot (U^{2m} \otimes U^{2n}).
		\end{equation}
	\end{itemize}
	Using these observations combined with Lemma \ref{lmma1}, we conclude that $\mathcal{B}=C(SU_q(2)\otimes SU_q(2)/\mathbb{Z}_2)$, which proves the claim.
\qed 

\blmma \label{lmma2}
Define the map $\lambda: \{z_i,y_i:1\leq i \leq 4\} \rightarrow \cll(\clh)$ given by 
\begin{IEEEeqnarray*}{rCll} 	z_1 & \mapsto &  a^*\otimes a^*, \qquad  & \qquad  y_1 \mapsto   b^* \otimes a^* \\ 	z_2 & \mapsto &  b \otimes a^*, \qquad & \qquad  y_2 \mapsto   a \otimes a^* \\ 	z_3 & \mapsto &  a^*\otimes b, \qquad  & \qquad  y_3 \mapsto   b \otimes b^* \\ 	z_4 & \mapsto &  b \otimes b, \qquad  & \qquad  y_4 \mapsto   a \otimes b.  \end{IEEEeqnarray*}.
Then the operators $\lambda(z_i)$ and $\lambda(y_i)$ satisfy the defining relations $\mathcal{R}_4$ of $C(SO_q(4))$.
\elmma 
\begin{proof}
The claim can be verified by a direct calculation. However, checking all the defining relations would be tedious. Instead, we use the following observations.
\begin{itemize}
	\item The operators $\Psi(z_i)$ and $\Psi(y_i)$ satisfy the defining relations $\mathcal{R}_4$.
		\item The
	$\ell^2(\mathbb{N}_0)\otimes \ell^2(\mathbb{N}_0)$-component  of $\Psi(z_i)$ and $\Psi(y_i)$  coincide with that of $\lambda(z_i)$ and $\lambda(y_i)$, respectively, for all $1\leq i\leq 4$.
	
	\item It therefore remains only to verify that, on
	$\ell^2(\mathbb{Z})\otimes \ell^2(\mathbb{Z})$, the operators appearing on the two sides of each relation in $\mathcal{R}_4$ agree for $\lambda(z_i)$ and $\lambda(y_i)$. This verification is straightforward, since no powers of $q$ are involved in this part of the computation.
	\end{itemize}
Hence, $\lambda(z_i)$ and $\lambda(y_i)$ satisfy all the defining relations $\mathcal{R}_4$ of $C(SO_q(4))$.
\end{proof}

It follows from Lemma \ref{lmma2}, that the map $\lambda$ extends uniquely to a surjective $C^*$-homomorphism 
$\Lambda: C(SO_q(4))\rightarrow C(SU_q(2)\otimes SU_q(2)/\mathbb{Z}_2)$.

\bthm \label{main} 	The $C^*$-homomorphism $\Lambda: C(SO_q(4))\rightarrow C((SU_q(2)\times SU_q(2))/\mathbb{Z}_2)$ is a quantum group isomorphism.
 \ethm
\begin{proof}
	We first show that $\Lambda$ is injective. By Proposition \ref{Hopf}, it suffices to show that the restriction 
$$
\Lambda:\operatorname{Pol}(SO_q(4))\longrightarrow
\mathcal{L}\bigl(L^2(\mathbb T)\otimes L^2(\mathbb T)\otimes
\ell^2(\mathbb N)\otimes\ell^2(\mathbb N)\bigr)
$$
is injective. 
	Let $P\in\operatorname{Pol}(SO_q(4))$ such that $\Lambda(P)=0$. Using the $q$-commutation relations among the generators, every monomial in $P$ can be arranged in the order
	$$
	z_1^{k_1}z_2^{k_2}z_3^{k_3}z_4^{k_4}
	y_1^{l_1}y_2^{l_2}y_3^{l_3}y_4^{l_4},
	\qquad k_i,l_i\in\mathbb{N}_0.
	$$
	Thus
	$$
	P=\sum_\lambda c_\lambda
	z_1^{k_{1,\lambda}}z_2^{k_{2,\lambda}}
	z_3^{k_{3,\lambda}}z_4^{k_{4,\lambda}}
	y_1^{l_{1,\lambda}}y_2^{l_{2,\lambda}}
	y_3^{l_{3,\lambda}}y_4^{l_{4,\lambda}}.
	$$
	By the explicit formulas for $\Lambda(z_i)$ and $\Lambda(y_i)$, we obtain
	$$
	\Lambda(P)=\sum_\lambda c_\lambda\,
	a^{n_\lambda}(a^*)^{m_\lambda}b^{r_\lambda}
	\otimes
	a^{n^\prime_\lambda}(a^*)^{m^\prime_\lambda}b^{r^\prime_\lambda}=0,
	$$
	where $n_{\lambda},m_\lambda,n^\prime_\lambda,m^\prime_\lambda\in\mathbb N_0$ and $r_\lambda,r^\prime_\lambda\in\mathbb Z$  such that $(n_{\lambda},m_{\lambda}, r_{\lambda}, n_{\lambda}^{\prime},m_{\lambda}^{\prime},r_{\lambda}^{\prime})\neq (n_{\mu},m_{\mu},r_{\mu},n_{\mu}^{\prime},m_{\mu}^{\prime},r_{\mu}^{\prime})$ if $\lambda\neq \mu$.
		Since
	$$
	\left\{
	a^n(a^*)^m b^k\otimes
	a^{n^\prime}(a^*)^{m^\prime}b^{k^\prime}:
	n,n^\prime,m,m^\prime\in\mathbb N_0,\ k,k^\prime\in\mathbb Z
	\right\}
	$$
	forms a linear basis of $\operatorname{Pol}(SU_q(2)\otimes SU_q(2))$, it follows that $c_\lambda=0$ for all $\lambda$. Hence $P=0$, and therefore $\ker\Lambda={0}$. Thus $\Lambda$ is injective.
	To prove the rest, we need to show that $\Lambda$ is a quantum group homomorphism. For that, we have 
	\begin{equation*}
		V = \begin{pmatrix}
			y_2^* & -q y_1^* & q^2 z_4^* & -q z_3^* \\
			-q y_4 & z_1^* & -q y_3 & -q z_2^* \\
			y_4 & y_3 & y_2 & y_1 \\
			z_4 & z_3 & z_2 & z_1
		\end{pmatrix}.
	\end{equation*}
Thus, 	we obtain the explicit component formulas:
	\begin{align*}
		\Delta(z_1) &= z_1 \otimes z_1 - q z_2 \otimes y_1 - q z_3 \otimes y_4^* + q^2 z_4 \otimes z_4^* \nonumber \\
		\Delta(z_2) &= z_1 \otimes z_2 + z_2 \otimes y_2 - q z_3 \otimes y_3 - q z_4 \otimes z_3^* \nonumber \\
		\Delta(z_3) &= z_1 \otimes z_3 + z_3 \otimes y_2^* - q z_2 \otimes y_3^* - q z_4 \otimes z_2^* \nonumber \\
		\Delta(z_4) &= z_1 \otimes z_4 + z_2 \otimes y_4 + z_3 \otimes y_1^* + z_4 \otimes z_1^* \nonumber \\
		\Delta(y_1) &= z_1 \otimes y_1 + y_1 \otimes y_2 - q z_3 \otimes y_3^* - q y_3^* \otimes y_4^* \nonumber \\
		\Delta(y_2) &= y_2 \otimes y_2 - q y_1 \otimes z_2 - q y_4 \otimes y_4^* + q^2 y_3^* \otimes y_3 \nonumber \\
		\Delta(y_3) &= y_2 \otimes y_3 + y_3 \otimes y_2^* + z_2 \otimes y_4^* + z_3^* \otimes y_1^* \nonumber \\
		\Delta(y_4) &= y_4 \otimes y_2 + z_1^* \otimes y_4 - q z_2^* \otimes z_4 - q y_3^* \otimes z_2
	\end{align*}
	One can  now  directly verify that $\Delta\circ \Lambda=(\Lambda\otimes \Lambda)\circ \Delta$, proving the claim.

\end{proof}

\section*{Acknowledgements} The author would like to thank Dr. Manabendra Giri for suggesting the reference \cite{Fio-1994aa} to obtain an explicit formula for the quantum determinant of $\mathcal{O}_q(4)$.

\bigskip

\noindent{\sc Bipul Saurabh} (\texttt{bipul.saurabh@iitgn.ac.in},  \texttt{saurabhbipul2@gmail.com})\\
{\footnotesize Department of Mathematics,  Indian Institute of Technology, Gandhinagar, Palaj, Gandhinagar 382055, India}

\end{document}